\documentclass[12pt,onecolumn,leqno]{article}

\usepackage{graphicx}
\usepackage{indentfirst}
\usepackage{amsmath}
\usepackage{amsfonts}
\usepackage{amssymb}
\usepackage{color}
\usepackage[mathscr]{eucal}
\usepackage{latexsym}

\newtheorem{theorem}{Theorem}[section]
\newtheorem{lemma}{Lemma}[section]

\newtheorem{prop}{Proposition}[section]

\numberwithin{equation}{section}

\title{$B(0,N)$-graded Lie superalgebras coordinatized by quantum tori}
    \author{Hongjia Chen,\quad  Yun Gao\footnote{Research was
partially supported by NSERC of Canada and Chinese Academy of
Science. This paper is dedicated to Professor Sheng Gong on the
occasion of his 75th birthday.}, \quad Shikui Shang}
    \date{}  

\begin{document}

\newcommand{\supercite}[1]{\textsuperscript{\cite{#1}}}

\maketitle              

     \begin{abstract}
          We use a fermionic extension of the bosonic module to obtain a class
          of $B(0,N)$-graded Lie superalgebras  with nontrivial central extensions.
     \end{abstract}


\setcounter{section}{-1}       

\section{Introduction}

$B(M-1, N)$-graded Lie superalgebras were first investigated and
classified up to central extension by Benkart-Elduque (see also
Garcia-Neher's work in [GN]). Those root graded Lie superalgebras
are a super-analog of root graded Lie algebras.
        Fermionic and bosonic representations
       for the affine Kac-Moody Lie algebras were studied by Frenkel
       [F1,2]
       and Kac-Peterson \cite{KP}. Feingold-Frenkel \cite{FF}
        constructed representations for all classical affine Lie algebras
       by using Clifford or Weyl algebras with infinitely many
       generators. They also obtained realizations for certain affine
       Lie superalgebras including the affine $B(0, N)$.
       \cite{G}
gave bosonic and fermionic representations for the extended affine
Lie
       algebra $\widetilde{gl_N(\mathbb{C}_q)}$, where $\mathbb{C}_q$ is the quantum
       torus in two variables. [CG] constructed  modules for some
       BC$_N$-graded Lie algebras by considering a fermionic
       extension of the fermionic module.

       In this paper, we will consider a fermionic extension of the bosonic
       module
       to obtain a class
          of $B(0,N)$-graded Lie superalgebras  with nontrivial central extensions.

       The organization of the paper is as follows. In Section 1, we review some basics
       on the quantum torus and present examples of $B(0,N)$-graded Lie superalgebras
       coordinatized by quantum tori
       which are subalgebras of $\widehat{gl(1,2N)}(\mathbb{C}_q)$. In Section 2, we use bosons
       and fermions to construct representations for those examples of $B(0,N)$-graded Lie superalgebras.

       Throughout this paper, we denote the field of complex numbers and the ring of integers by
       $\mathbb{C}$ and $\mathbb{Z}$ respectively.

\section{B(0,N)-graded Lie superalgebras}

       We first recall some basics on quantum
       tori and  then
       go on to present examples of B(0,N)-graded Lie
       superalgebras coordinatized by quantum tori. For more information on Lie superalgebras
       graded by root systems, see \cite{BE1}-\cite{BE2} and [GN].

       Let $q$ be a non-zero complex number. A quantum torus
       associated to $q$ (see \cite{M}) is the unital  associative
       $\mathbb{C}$-algebra $\mathbb{C}_q[x^{\pm},y^{\pm}]$(or,
       simply $\mathbb{C}_q$) with generators $x^{\pm}$, $y^{\pm}$
       and relations
           \begin{equation}
                xx^{-1}=x^{-1}x=yy^{-1}=y^{-1}y=1 \quad and \quad yx=qxy.
           \end{equation}
       Set $\Lambda(q)=\{n \in \mathbb{Z}|q^n=1\}$. From \cite{BGK}
       we see that $[\mathbb{C}_q,\mathbb{C}_q]$ has a basis consisting
       of monomials $x^my^n$ for $m \notin \Lambda(q)$ or $n \notin
       \Lambda(q)$.

       Let $\mbox{ }\bar{} \mbox{ }$ be the anti-involution on
       $\mathbb{C}_q$ given by
           \begin{equation}
                \bar{x}=x,\quad \bar{y}=y^{-1}.
           \end{equation}
       We have $\mathbb{C}_q=\mathbb{C}_q^+\oplus\mathbb{C}_q^-$,
       where $\mathbb{C}_q^{\pm}=\{s \in \mathbb{C}_q|\bar{s}=\pm
       s\}$, then    \\
       \parbox{1cm}{\begin{eqnarray}\end{eqnarray}}
                            \hfill \parbox{13.66cm}
              {\begin{eqnarray*}
                    & \mathbb{C}_q^+=span\{x^my^n +
                      \overline{x^my^n}|m \in \mathbb{Z},n \geq 0\}, & \\
                    & \mathbb{C}_q^-=span\{x^my^n -
                      \overline{x^my^n}|m \in \mathbb{Z},n>0\}. & \hfill
                \end{eqnarray*}}

       Let $M,N$ be two positive integers. We have a Lie
       superalgebra $gl(M,N)(\mathbb{C}_q)$ of $(M+N)$ by $(M+N)$
       matrices with entries from $\mathbb{C}_q$.

       We form a central extension of Lie superalgebra $gl(M,N)(\mathbb{C}_q)$
       as was done in \cite{G} and [CG].
           \begin{equation}
               \widehat{gl(M,N)}(\mathbb{C}_q)=gl(M,N)(\mathbb{C}_q)
               \oplus \Bigl(\sum_{n \in \Lambda(q)}\oplus
               \mathbb{C}c(n) \Bigr) \oplus \mathbb{C}c_y
           \end{equation}
       with bracket \\
         \parbox{1cm}{\begin{eqnarray}\end{eqnarray}}
                                      \hfill \parbox{13.66cm}
              {\begin{eqnarray*}
                     [A(x^my^n),B(x^py^s)]_s=A(x^my^n)B(x^py^s)-(-1)^{degAdegB}B(x^py^s)A(x^my^n)&& \\
                           +m q^{np} str(AB) \delta_{m+p,0}\delta_{\overline{n+s},\overline{0}} c(n+s)
                                  +n q^{np}str(AB) \delta_{m+p,0}\delta_{n+s,0} c_y && \hfill
               \end{eqnarray*}}
       for $m,p,n,s \in \mathbb{Z},$ $A,B \in gl(M,N)_{\alpha},\alpha=
       \bar{0} \mbox{ or } \bar{1}$, where $str$ is the super-trace
       of the Lie superalgebra $gl(M, N)$,
       $c(u)$ with $u \in \Lambda(q)$ and $c_y$ are central
       elements of $\widehat{gl(M,N)}(\mathbb{C}_q)$, $\overline{t} \mbox{ means }
       \overline{t} \in \mathbb{Z}/\Lambda(q),\mbox{for }t \in \mathbb{Z}$.

       Now we present the examples of Lie superalgebra graded by the
       root system of type $B(0,N)$.  We first set
           \begin{equation*}J=\left(
                 \begin{array}{ccc}
                      1 & 0  \\
                      0 & -I_{2N}
                 \end{array}
                      \right),
                      G=\left(
                 \begin{array}{ccc}
                      1 & 0 & 0 \\
                      0 & 0 & I_N \\
                      0 & -I_N & 0
                 \end{array}
                      \right) \in M_{2N+1}(\mathbb{C}_q).
           \end{equation*}
       Then, $G$ and $J$ are invertible {\small $ (2N+1) \times
       (2N+1)$}-matrices. Using the matrix $G$ and $J$, we define a
       superspace $\mathcal{S}$ with:
           \begin{equation*}
                \mathcal{S}_{\bar{0}}=\{X \in gl(1,2N)
                   (\mathbb{C}_q)_{\bar{0}}|\bar{X}^tG+GX=0\}
           \end{equation*}
           \begin{equation*}
                \mathcal{S}_{\bar{1}}=\{X \in gl(1,2N)
                   (\mathbb{C}_q)_{\bar{1}}|\bar{X}^tG-JGX=0\}
           \end{equation*}

       We can easily see that $\mathcal{S}$ is a subalgebra of
       $gl(1,2N)(\mathbb{C}_q)$ over $\mathbb{C}$. The general form of
       a matrix in $\mathcal{S}$ is
           \begin{equation}
                 \left(
           \begin{array}{ccc}
                 a & b & c \\
                 \bar{c}^t  & A & S \\
                 -\bar{b}^t & T & -\bar{A}^t\\
            \end{array}
                 \right) \quad \mbox{with  }\bar{a}=-a \quad
                 \bar{S}^t=S \quad \mbox{and} \quad \bar{T}^t=T,
            \end{equation}
       where $A,S,T$ are $N\times N$ sub-matrices. Then the Lie superalgebra
       $\mathcal{G}=[\mathcal{S},\mathcal{S}]_s$, is a B(0,N)-graded
       Lie superalgebra.

       As in \cite{AABGP}, we easily
       know that:
            $$\mathcal{G}=\{Y \in
                 gl(1,2N)(\mathbb{C}_q)|str(Y)\equiv 0 \mbox{ mod }
                 [\mathbb{C}_q,\mathbb{C}_q]\}.$$

       Putting
            \begin{equation}
                 \mathcal{H}=\Bigl\{\sum_{i=1}^Na_i(e_{ii}
                  -e_{N+i,N+i})|a_i \in \mathbb{C}\Bigr\},
            \end{equation}
       then $\mathcal{H}$ is a $N$-dimensional abelian subalgebra
       of $\mathcal{G}$. Defining $\delta_i \in
       \mathcal{H}^*,i=1,\cdots,N$, by
            \begin{equation}
                 \delta_i\biggl(\sum_{j=1}^Na_j(e_{jj}-e_{N+j,N+j})\biggr)=a_i
            \end{equation}
       for $i=1,\cdots,N$. Setting $\mathcal{G}_{\alpha}=\{x\in
       \mathcal{G}|[h,x]_s=\alpha(h)x,\mbox{ for all } h \in
       \mathcal{H}\}$ as usual, we have
            \begin{equation}
                 \mathcal{G}=\mathcal{G}_0 \oplus \sum_{i \neq j}
                 \mathcal{G}_{\delta_i-\delta_j}
                 \oplus \sum_{i<j}(\mathcal{G}_{\delta_i+\delta_j}
                 \oplus \mathcal{G}_{-\delta_i-\delta_j})
                 \oplus \sum_i (\mathcal{G}_{\delta_i}
                 \oplus \mathcal{G}_{-\delta_i}
                 \oplus \mathcal{G}_{2\delta_i}
                 \oplus \mathcal{G}_{-2\delta_i}),
            \end{equation}
       where \\
          \parbox{1cm}{\begin{eqnarray}\end{eqnarray}}\hfill \parbox{13.66cm}
             {\begin{eqnarray*}
               &\mathcal{G}_{\delta_i-\delta_j}=span_{\mathbb{C}}
                 \{\tilde{f}_{ij}(m,n)=x^my^ne_{ij}-
                 \overline{x^my^n}e_{N+j,N+i} |m,n \in \mathbb{Z}\}, & \\
               &\mathcal{G}_{\delta_i+\delta_j}=span_{\mathbb{C}}\{
                 \tilde{g}_{ij}(m,n)=x^my^ne_{i,N+j}+
                 \overline{x^my^n}e_{j,N+i}|m,n \in \mathbb{Z} \},& \\
               &\mathcal{G}_{-\delta_i-\delta_j}=span_{\mathbb{C}}
                 \{ \tilde{h}_{ij}(m,n)=-x^my^ne_{N+i,j}-
                 \overline{x^my^n}e_{N+j,i} |m,n \in \mathbb{Z} \}, & \\
               &\mathcal{G}_{2\delta_i}=span_{\mathbb{C}}\{\tilde{g}_{ii}(m,n)=
                 (x^my^n+\overline{x^my^n})e_{i,N+i}|m,n \in \mathbb{Z} \}, & \\
               &\mathcal{G}_{-2\delta_i}=span_{\mathbb{C}}\{\tilde{h}_{ii}(m,n)=
                 -(x^my^n+\overline{x^my^n})e_{N+i,i}|m,n \in \mathbb{Z} \},&\\
               &\mathcal{G}_{\delta_i}=span_{\mathbb{C}}\{\tilde{e}_i(m,n)=
                 -x^my^ne_{i,0}-\overline{x^my^n}e_{0,N+i}|m,n \in \mathbb{Z} \},& \\
               &\mathcal{G}_{-\delta_i}=span_{\mathbb{C}}\{\tilde{e}_i^*(m,n)=
                 x^my^ne_{N+i,0}-\overline{x^my^n}e_{0,i}|m,n \in \mathbb{Z}\} &\hfill
             \end{eqnarray*}}
       and
            \begin{equation*}
                 \mathcal{G}_0=span_{\mathbb{C}}\{\tilde{f}_{ii}(m,n)
                 +\tilde{e}_0(m,n),\tilde{e}_0(p,s)|1 \leq i \leq N,m,n
                 \in \mathbb{Z},p \notin \Lambda(q)\mbox{ or } s \notin \Lambda(q)\}
            \end{equation*}
       where $\tilde{e}_0(m,n)=-(x^my^n-\overline{x^my^n})e_{0,0}$.

       We then have a central extension of $\mathcal{G}$
            \begin{equation}
                 \widehat{\mathcal{G}}=\mathcal{G} \oplus
                 \Bigl(\sum_{n \in \Lambda(q)}\oplus
                 \mathbb{C}c(n) \Bigr) \oplus \mathbb{C}c_y
            \end{equation}
       with bracket as (1.5).

       We have
            \begin{prop}
                 \begin{equation}
                      [\tilde{g}_{ij}(m,n),\tilde{g}_{kl}(p,s)]_s=0
                 \end{equation}
                 \begin{equation}
                      [\tilde{g}_{ij}(m,n),\tilde{f}_{kl}(p,s)]_s=
                      -\delta_{il}q^{ms}\tilde{g}_{kj}(m+p,n+s)
                      -\delta_{jl}q^{(s-n)m}\tilde{g}_{ki}(m+p,s-n)
                 \end{equation}
            \begin{eqnarray}
                 &&[\tilde{g}_{ij}(m,n),\tilde{h}_{kl}(p,s)]_s \nonumber \\
                 &=&-\delta_{ik}q^{-n(m+p)}\tilde{f}_{jl}(m+p,s-n)
                 -\delta_{jk}q^{np}\tilde{f}_{il}(m+p,n+s)  \nonumber  \\
                 &&-\delta_{il}q^{-(mn+np+ps)}\tilde{f}_{jk}(m+p,-(n+s))
                 -\delta_{jl}q^{(n-s)p}\tilde{f}_{ik}(m+p,n-s)      \\
                 &&+m q^{np}\delta_{jk}\delta_{il}\delta_{m+p,0}\delta_{\overline{n+s},\overline{0}}(c(n+s)+c(-n-s))
                         \nonumber \\
                 &&+m\delta_{ik}\delta_{jl}\delta_{m+p,0}\delta_{\overline{n-s},\overline{0}}(c(n-s)+c(s-n)) \nonumber
            \end{eqnarray}
            \begin{equation}
                 [\tilde{g}_{ij}(m,n),\tilde{e}_{k}(p,s)]_s=[\tilde{g}_{ij}(m,n),\tilde{e}_{0}(p,s)]_s=0
            \end{equation}
            \begin{equation}
                 [\tilde{g}_{ij}(m,n),\tilde{e}_k^*(p,s)]_s=-\delta_{ik}q^{-n(m+p)}
                 \tilde{e}_j(m+p,s-n)-\delta_{jk}q^{np}\tilde{e}_i(m+p,n+s)
            \end{equation}
            \begin{equation}
                 [\tilde{g}_{ij}(m,n),\tilde{e}_{0}(p,s)]_s=0
            \end{equation}
            \begin{eqnarray}
                 [\tilde{f}_{ij}(m,n),\tilde{f}_{kl}(p,s)]_s&=&\delta_{jk}q^{np}
                 \tilde{f}_{il}(m+p,n+s)-\delta_{il}q^{sm}\tilde{f}_{kj}(m+p,n+s) \nonumber  \\
                  &&-2mq^{np}\delta_{jk}\delta_{il}\delta_{m+p,0}\delta_{\overline{n+s},\overline{0}}c(n+s)
            \end{eqnarray}
            \begin{equation}
                 [\tilde{f}_{ij}(m,n),\tilde{h}_{kl}(p,s)]_s=-\delta_{ik}q^{-n(m+p)}
                 \tilde{h}_{jl}(m+p,s-n)-\delta_{il}q^{ms}\tilde{h}_{kj}(m+p,n+s)
            \end{equation}
            \begin{equation}
                 [\tilde{f}_{ij}(m,n),\tilde{e}_k(p,s)]_s=\delta_{jk}q^{np}\tilde{e}_i(m+p,n+s)
            \end{equation}
            \begin{eqnarray}
                 &[\tilde{f}_{ij}(m,n),\tilde{e}_k^*(p,s)]_s
                 =-\delta_{ik}q^{-n(m+p)}\tilde{e}_j^*(m+p,s-n)&
            \end{eqnarray}
            \begin{equation}
                 [\tilde{f}_{ij}(m,n),\tilde{e}_{0}(p,s)]_s=0
            \end{equation}
            \begin{equation}
                 [\tilde{h}_{ij}(m,n),\tilde{h}_{kl}(p,s)]_s=0
            \end{equation}
            \begin{equation}
                 [\tilde{h}_{ij}(m,n),\tilde{e}_k(p,s)]_s=\delta_{jk}q^{np}\tilde{e}_i^*(m+p,n+s)+
                   \delta_{ik}q^{-n(m+p)}\tilde{e}_j^*(m+p,s-n)
            \end{equation}
            \begin{equation}
                 [\tilde{h}_{ij}(m,n),\tilde{e}^*_{k}(p,s)]_s=0
            \end{equation}
            \begin{equation}
                 [\tilde{h}_{ij}(m,n),\tilde{e}_{0}(p,s)]_s=0
            \end{equation}
            \begin{equation}
                 [\tilde{e}_i(m,n),\tilde{e}_k(p,s)]_s=q^{m(s-n)}\tilde{g}_{ki}(m+p,s-n)
            \end{equation}
            \begin{eqnarray}
                 &&[\tilde{e}_i(m,n),\tilde{e}_k^*(p,s)]_s \nonumber  \\
                 &=&\delta_{ik}q^{-n(m+p)}\tilde{e}_0(m+p,s-n)+q^{p(n-s)}\tilde{f}_{ik}(m+p,n-s) \\
                 &&-m\delta_{ik}\delta_{m+p,0}\delta_{\overline{n-s},\overline{0}}(c(n-s)+c(s-n)) \nonumber
            \end{eqnarray}
            \begin{equation}
                 [\tilde{e}_i(m,n),\tilde{e}_0(p,s)]_s=-q^{np}\tilde{e}_i(m+p,n+s)+q^{p(n-s)}\tilde{e}_i(m+p,n-s)
            \end{equation}
            \begin{equation}
                 [\tilde{e}_i^*(m,n),\tilde{e}_k^*(p,s)]_s=q^{m(s-n)}\tilde{h}_{ki}(m+p,s-n)
            \end{equation}
            \begin{equation}
                 [\tilde{e}_i^*(m,n),\tilde{e}_0(p,s)]_s=-q^{np}\tilde{e}_i^*(m+p,n+s)+q^{p(n-s)}\tilde{e}_i^*(m+p,n-s)
            \end{equation}
            \begin{eqnarray}
                 &&[\tilde{e}_0(m,n),\tilde{e}_0(p,s)]_s  \\
                 &=&-(q^{np}-q^{sm})\tilde{e}_0(m+p,n+s)-(q^{m(s-n)}-q^{-n(m+p)})\tilde{e}_0(m+p,s-n) \nonumber \\
                 &&+mq^{np}\delta_{m+p,0}\delta_{\overline{n+s},\overline{0}}(c(n+s)+c(-n-s))\nonumber  \\
                 &&-m\delta_{m+p,0}\delta_{\overline{n-s},\overline{0}}(c(n-s)+c(s-n)) \nonumber
            \end{eqnarray}
       for all $m,p,n,s \in \mathbb{Z}$ and $1 \leq i,j,k,l \leq N$.
            \end{prop}

            \noindent{\bf Remark 1.1} For our situation in the decomposition of $B(0, N)$-graded Lie superalgebras
            in \S 5[Be2], $A=\mathbb{C}_q^+$ and $B=\mathbb{C}_q^-$.

\section{The module construction}

       In this section, we follow the method in
       \cite{G} and
       \cite{CG} to construct representations for the Lie
       superalgebras which are given in Section 1. The idea goes
       back to \cite{FF}.

       Let $\mathcal{R}$ be an associative algebra. Let $\rho=\pm
       1$. We define a $\rho$-bracket on $\mathcal{R}$ as follow:
            \begin{equation}
                 \{a,b\}_{\rho}=ab+\rho ba,\quad a,b \in \Re.
            \end{equation}
       It is easy to see that
            \begin{equation}
                 \{a,b\}_{\rho}=\rho \{b,a\}_{\rho}
            \end{equation}
       for $a,b,c \in \mathcal{R}$.

       Define $\mathfrak{a}$ to be the unital associative algebra
       with $2N$ generators $a_i,a_i^*,1 \leq i \leq N $, subject to
       relations
            \begin{equation}
                 \{a_i,a_j\}_{-}=\{a_i^*,a_j^*\}_{-}=0, \quad \mbox{and}
                                    \quad \{a_i,a_j^*\}_{-}=-\delta_{ij}.
            \end{equation}
       Let the associative algebra $\alpha (N,-)$ be generated by
            \begin{equation}
                 \{u(m)|u \in \bigoplus_{i=1}^N (\mathbb{C}a_i
                 \oplus \mathbb{C}a^*_i),m \in \mathbb{Z}\}
            \end{equation}
       with the relations
            \begin{equation}
                 \{u(m),v(n)\}_{-}=\{u,v\}_{-}\delta_{m+n,0}.
            \end{equation}
       The normal ordering is defined as in \cite{FF}(see also \cite{F2} or \cite{CG}). \\
            \parbox{1cm}{\begin{eqnarray}\end{eqnarray}}\hfill \parbox{13.66cm}
               {\begin{eqnarray*} :u(m)v(n):&=&\left\{
                   \begin{array}
                        {r@{ }l}  & u(m)v(n) \hspace{3.22cm}                        \mbox{if } n>m,
                               \\ & \frac{1}{2}\bigl(u(m)v(n)+v(n)u(m)\bigr)  \quad \mbox{if } n=m,
                               \\ &v(n)u(m) \hspace{3.25cm}                         \mbox{if } n<m,
                   \end{array} \right.  \\&=&:v(n)u(m): \hfill
                \end{eqnarray*}}
       for $n,m \in \mathbb{Z},u,v \in \mathfrak{a}.$

       Next we consider an extension of the algebra $\alpha (N,-)$. The generators
            \begin{equation}
                 \{e(m)|m \in \mathbb{Z}\}
            \end{equation}
       span an infinite-dimensional Clifford algebra with relations
            \begin{equation}
                 \{e(m),e(n)\}_{+}=e(m)e(n)+e(n)e(m)=-\delta_{n+m,0}.
            \end{equation}
       Let $\alpha_\tau(N)$ denote the algebra obtained by adjoining to
       $\alpha (N,-)$ the generators (2.7) with relations (2.8) and
            \begin{equation}
                 \{a_i(m),e(n)\}_{\tau}=0=\{a_i^*(m),e(n)\}_{\tau},
                 \mbox{  for } \tau=\pm1.
            \end{equation}

       The normal ordering is given as follows \\
            \parbox{1cm}{\begin{eqnarray}\end{eqnarray}}\hfill \parbox{13.66cm}
               {\begin{eqnarray*}& :e(m)e(n):=\left\{
                  \begin{array}
                        {r@{ }l}
                        & e(m)e(n) \hspace{3.18cm} \mbox{if }        n>m,\vspace{0.1cm}
                        \\ &\frac{1}{2}\bigl(e(m)e(n)-e(n)e(m)\bigr) \quad \mbox{if }  n=m,
                        \\& -e(n)e(m) \hspace{2.88cm} \mbox{if }n<m,
                  \end{array} \right., & \\
                            & :a_i(m)e(n):=a_i(m)e(n)=-\tau e(n)a_i(m),& \\
                            & :a_i^*(m)e(n):=a_i^*(m)e(n)=-\tau e(n)a_i^*(m),&  \hfill
                \end{eqnarray*}}
       for $n,m \in \mathbb{Z},1 \leq i,j \leq N.$ Set
            \begin{equation}
                 \theta(n)=\left\{
                   \begin{array}
                        {r@{ \quad }l}
                         1, & \mbox{for } n>0, \\
                         \frac{1}{2} ,& \mbox{for } n=0,  \\
                         0, & \mbox{for } n<0,
                   \end{array} \right. \mbox{   then  }     1-\theta(n)=\theta(-n).
            \end{equation}
       We have \\
            \parbox{1cm}{\begin{eqnarray}\end{eqnarray}}\hfill \parbox{13.66cm}
               {\begin{eqnarray*}& :a_i(m)a_j(n):=a_i(m)a_j(n)= a_j(n)a_i(m),& \\
                                  & :a_i^*(m)a_j^*(n):=a_i^*(m)a_j^*(n)=a_j^*(n)a_i^*(m),& \hfill
                \end{eqnarray*}}
       and \\
            \parbox{1cm}{\begin{eqnarray}\end{eqnarray}}\hfill \parbox{13.66cm}
               {\begin{eqnarray*}
                     &a_i(m)a_j^*(n)=:a_i(m)a_j^*(n):-\delta_{ij}\delta_{m+n,0}   \theta(m-n), & \\
                     & a_j^*(n)a_i(m)=:a_i(m)a_j^*(n):+\delta_{ij}\delta_{m+n,0}\theta(n-m),& \\
                     & e(m)e(n)=:e(m)e(n):-\delta_{m+n,0} \theta(m-n). &
                                  \hfill
                \end{eqnarray*}}

       Let $\alpha(N,-)^+$ be the subalgebra generated by
       $a_i(n),a_j^*(m),a_k^*(0)$, for $n,m>0$, and $1 \leq i,j,k
       \leq N.$ Let $\alpha(N,-)^-$ be the subalgebra generated by
       $a_i(n),a_j^*(m),a_k(0),$ for $n,m<0,$ and $1 \leq i,j,k \leq
       N.$ Those generators in $\alpha(N,-)^+$ are called
       annihilation operators while those in $\alpha(N,-)^-$ are
       called creation operators. Let $V(N,-)$ be a simple
       $\alpha(N,-)$-module containing an element $v_0$, called a
       ``vacuum vector$"$, and satisfying
            \begin{equation}
                 \alpha(N,-)^+v_0=0.
            \end{equation}
       So all annihilation operators kill $v_0$ and
            \begin{equation}
                 V(N,-)=\alpha(N,-)^-v_0.
            \end{equation}

       Let $V_0$ be a simple Clifford module for the Clifford
       algebra generated by (2.7) with relations (2.8) and
       containing ``vacuum vector" $v'_0$, which is killed by
       annihilation operators.(Here we call $e(m)$ annihilation
       operator if $m>0$, or a creation operator if $m<0$. e(0)
       acts as scalar.) Because of (2.9), we see that the
       $\alpha_\tau(N)$-module
            \begin{equation}
                 V_\tau(N)=V(N,-) \otimes V_0=\alpha_\tau(N)v'_0
            \end{equation}
       is simple.

       Now we define our operators on $V_\tau(N)$. For any
       $m,n \in \mathbb{Z},1 \leq i,j \leq N,$ set

\begin{eqnarray}
                 & f_{ij}(m,n)=\sum_{s \in \mathbb{Z}}
                 q^{-ns}:a_i(m-s)a_j^*(s):\\
                 & g_{ij}(m,n)=\sum_{s \in \mathbb{Z}}
                 q^{-ns}:a_i(m-s)a_j(s):\\
                 & h_{ij}(m,n)=\sum_{s \in \mathbb{Z}}
                 q^{-ns}:a_i^*(m-s)a_j^*(s):\\
                & e_i(m,n)=\sum_{s \in \mathbb{Z}}
                q^{-ns}:a_i(m-s)e(s):\\
                & e_i^*(m,n)=\sum_{s \in \mathbb{Z}}
                q^{-ns}:a_i^*(m-s)e(s):\\
                 & e_0(m,n)=\sum_{s \in \mathbb{Z}} q^{-ns}:e(m-s)e(s):.
\end{eqnarray}

We will list only those commutation relations involving $e(s)$ (See
[CG]).
            \begin{lemma} We have

                 \begin{equation}
                      [a_i(m)a_j(n),a_k(p)e(s)]=[a_i(m)a_j(n),e(p)e(s)]=0
                 \end{equation}
                 \begin{equation}
                      [a_i(m)a_j(n),a^*_k(p)e(s)]
                      =-\delta_{ik}\delta_{m+p,0}a_j(n)e(s)-\delta_{jk}\delta_{n+p,0}a_i(m)e(s)
                 \end{equation}

                 \begin{equation}
                      [a_i(m)a^*_j(n),a_k(p)e(s)]=\delta_{jk}\delta_{n+p,0}a_i(m)e(s),
                 \end{equation}
                 \begin{equation}
                      [a_i(m)a^*_j(n),a_k^*(p)e(s)]=-\delta_{ik}\delta_{m+p,0}a_j^*(n)e(s),
                 \end{equation}
                 \begin{equation}
                      [a_i(m)a^*_j(n),e(p)e(s)]=[a^*_i(m)a^*_j(n),a^*_k(p)a^*_l(s)]=0
                 \end{equation}
                 \begin{equation}
                      [a^*_i(m)a^*_j(n),a_k(p)e(s)]
                      =\delta_{jk}\delta_{n+p,0}a_i^*(m)e(s)+\delta_{ik}\delta_{m+p,0}a_j^*(m)e(s),
                 \end{equation}
                 \begin{equation}
                      [a^*_i(m)a^*_j(n),a^*_k(p)e(s)]=[a^*_i(m)a^*_j(n),e(p)e(s)]=0
                 \end{equation}
                 \begin{equation}
                      \{a_i(m)e(n),a_k(p)e(s)\}_+=\tau\delta_{n+s,0}a_i(m)a_k(p)
                 \end{equation}
                 \begin{equation}
                      \{a_i(m)e(n),a_k^*(p)e(s)\}_+
                      =\tau\delta_{n+s,0}a_k^*(p)a_i(m)+\tau\delta_{ik}\delta_{m+p,0}e(n)e(s)
                 \end{equation}
                 \begin{equation}
                      [a_i(m)e(n),e(p)e(s)]=\delta_{n+s,0}a_i(m)e(p)-\delta_{n+p,0}a_i(m)e(s)
                 \end{equation}
                 \begin{equation}
                      \{a_i^*(m)e(n),a_k^*(p)e(s)\}_+=\tau\delta_{n+s,0}a_i^*(m)a_k^*(p)
                 \end{equation}
                 \begin{equation}
                      [a_i^*(m)e(n),e(p)e(s)]=\delta_{n+s,0}a_i^*(m)e(p)-\delta_{n+p,0}a_i^*(m)e(s)
                 \end{equation}
                 \begin{eqnarray}
                      &&[e(m)e(n),e(p)e(s)]    \\
                      &=&-\delta_{n+p,0}e(m)e(s)+\delta_{m+p,0}e(n)e(s)-\delta_{n+s,0}e(p)e(m)+\delta_{m+s,0}e(p)e(n)
                      \nonumber
                 \end{eqnarray}
       for m,n,p,s $\in \mathbb{Z}$ and $1 \leq i,j,k,l \leq N$.
            \end{lemma}
       \textbf{$Proof$ } We only check some of them.
\begin{eqnarray*}
                 &&\{a_i(m)e(n),a_k(p)e(s)\}_+
                    =a_i(m)e(n)a_k(p)e(s)+a_k(p)e(s)a_i(m)e(n) \qquad   \\
                 &=&-\tau\bigl(a_i(m)a_k(p)e(n)e(s)+a_k(p)a_i(m)e(s)e(n)\bigr) \\
                 &=&-\tau\bigl(a_i(m)a_k(p)(e(n)e(s)+e(s)e(n))\bigr)              \\
                 &=&\tau\delta_{n+s,0}a_i(m)a_k(p);
            \end{eqnarray*}
            \begin{eqnarray*}
                 &&\{a_i(m)e(n),a_k^*(p)e(s)\}_+
                    =a_i(m)e(n)a_k^*(p)e(s)+a_k^*(p)e(s)a_i(m)e(n)   \qquad \qquad \\
                 &=&-\tau\bigl(a_i(m)a_k^*(p)e(n)e(s)+a_k^*(p)a_i(m)e(s)e(n)\bigr)     \\
                 &=&-\tau\bigl(a_k^*(p)a_i(m)-\delta_{ik}\delta_{m+p,0}\bigr)e(n)e(s)-\tau a_k^*(p)a_i(m)e(s)e(n) \\
                 &=&\tau\delta_{n+s,0}a_k^*(p)a_i(m)
                 +\tau\delta_{ik}\delta_{m+p,0}e(n)e(s);
            \end{eqnarray*}
            \begin{eqnarray*}
                 &&[a_i(m)e(n),e(p)e(s)]=a_i(m)e(n)e(p)e(s)-e(p)e(s)a_i(m)e(n)    \qquad  \qquad   \\
                 &=&a_i(m)\bigl(e(n)e(p)e(s)-(-\tau)^2e(p)e(s)e(n)\bigr)                          \\
                 &=&a_i(m)\bigl(-\delta_{n+p,0}e(s)-e(p)e(n)e(s)-e(p)e(s)e(n)\bigr)                          \\
                 &=&a_i(m)\bigl(\delta_{n+s,0}e(p)-\delta_{n+p,0}e(s)\bigr)                                       \\
                 &=&\delta_{n+s,0}a_i(m)e(p)-\delta_{n+p,0}a_i(m)e(s);
            \end{eqnarray*}
            \begin{eqnarray*}
                 &&[e(m)e(n),e(p)e(s)]=e(m)e(n)e(p)e(s)-e(p)e(s)e(m)e(n)     \\
                 &=&e(m)(-\delta_{n+p,0}-e(p)e(n))e(s)-e(p)e(s)e(m)e(n)        \\
                 &=&-\delta_{n+p,0}e(m)e(s)-e(m)e(p)e(n)e(s)-e(p)e(s)e(m)e(n)             \\
                 &=&-\delta_{n+p,0}e(m)e(s)-(-\delta_{m+p,0}-e(p)e(m))e(n)e(s)-e(p)e(s)e(m)e(n)       \\
                 &=&-\delta_{n+p,0}e(m)e(s)+\delta_{m+p,0}e(n)e(s)+e(p)\bigl(e(m)e(n)e(s)-e(s)e(m)e(n)\bigr)  \\
                 &=&-\delta_{n+p,0}e(m)e(s)+\delta_{m+p,0}e(n)e(s)
                                   +e(p)\bigl(\delta_{m+s,0}e(n)-\delta_{n+s,0}e(m)\bigr)\\
                 &=&-\delta_{n+p,0}e(m)e(s)+\delta_{m+p,0}e(n)e(s)-\delta_{n+s,0}e(p)e(m)+\delta_{m+s,0}e(p)e(n).
            \end{eqnarray*}
       So (2.30),(2.31),(2.32) and (2.35) hold true.  $\hfill \blacksquare$

       In what follows we shall mean $\frac{q^{mn}-1}{q^n-1}=m$ if $n \in \Lambda(q)$.
       This will make our formula more concise.

       Next we list all brackets that are needed.  For all $m,p,n,s \in \mathbb{Z}$
       and $1 \leq i,j,k,l \leq N$,
            \begin{prop}
                 \begin{equation*}
                      [g_{ij}(m,n),g_{kl}(p,s)]=0
                 \end{equation*}

            \end{prop}
            \begin{prop}
                 \begin{equation*}
                      [g_{ij}(m,n),f_{kl}(p,s)]=-\delta_{il}q^{ms}g_{kj}(m+p,n+s)- \delta_{jl}q^{(s-n)m}g_{ki}(m+p,s-n)
                 \end{equation*}
            \end{prop}
            \begin{prop}
               \begin{eqnarray*}
                  &&[g_{ij}(m,n),h_{kl}(p,s)]\\
                  &=&-\delta_{ik}q^{-n(m+p)}f_{jl}(m+p,s-n)-\delta_{jk}q^{np}f_{il}(m+p,n+s)  \\
                  &&-\delta_{il}q^{-(mn+np+ps)}f_{jk}(m+p,-(n+s))-\delta_{jl}q^{(n-s)p}f_{ik}(m+p,n-s) \\
                  &&+\delta_{ik}\delta_{jl}\delta_{m+p,0}\frac{1}{2}(q^{s-n}+1)\frac{q^{m(s-n)}-1}{q^{s-n}-1}
                    +\delta_{jk}\delta_{il}\delta_{m+p,0}q^{np}\frac{1}{2}(q^{s+n}+1)\frac{q^{m(s+n)}-1}{q^{s+n}-1}
               \end{eqnarray*}
            \end{prop}
            \begin{prop}
                 \begin{equation*}
                      [g_{ij}(m,n),e_k(p,s)]=[g_{ij}(m,n),e_0(p,s)]=0
                 \end{equation*}
                 \begin{equation*}
                      [g_{ij}(m,n),e_k^*(p,s)]=-\delta_{ik}q^{-n(m+p)}e_j(m+p,s-n)-\delta_{jk}q^{np}e_i(m+p,n+s)
                 \end{equation*}
            \end{prop}
            \begin{prop}
                 \begin{eqnarray*}
                      \hspace{-0.3cm}[f_{ij}(m,n),f_{kl}(p,s)]
                       &=&\delta_{jk}q^{np}f_{il}(m+p,n+s)-\delta_{il}q^{sm}f_{kj}(m+p,n+s)\\
                       &&-\delta_{jk} \delta_{il} q^{np} \delta_{m+p,0} \frac{1}{2}
                         (q^{s+n}+1)\frac{q^{m(s+n)}-1}{q^{s+n}-1}
                 \end{eqnarray*}
            \end{prop}
            \begin{prop}
               \begin{eqnarray*}
                  &[f_{ij}(m,n),h_{kl}(p,s)]=-\delta_{ik}q^{-n(m+p)}h_{jl}(m+p,s-n)-\delta_{il}q^{ms}h_{kj}(m+p,n+s)&
               \end{eqnarray*}
            \end{prop}
            \begin{prop}
                 \begin{equation*}
                      [f_{ij}(m,n),e_k(p,s)]=\delta_{jk}q^{np}e_i(m+p,n+s)
                 \end{equation*}
                 \begin{equation*}
                      [f_{ij}(m,n),e_k^*(p,s)]=-\delta_{ik}q^{-n(m+p)}e_j^*(m+p,s-n)
                 \end{equation*}
                 \begin{equation*}
                      [f_{ij}(m,n),e_0(p,s)]=0
                 \end{equation*}
            \end{prop}
            \begin{prop}
                 \begin{equation*}
                      [h_{ij}(m,n),h_{kl}(p,s)]=0
                 \end{equation*}
            \end{prop}
            \begin{prop}
                 \begin{equation*}
                      [h_{ij}(m,n),e_k(p,s)]=\delta_{jk}q^{np}e_i^*(m+p,n+s)+\delta_{ik}q^{-n(m+p)}e_j^*(m+p,s-n)
                 \end{equation*}
                 \begin{equation*}
                      [h_{ij}(m,n),e_k^*(p,s)]=[h_{ij}(m,n),e_0(p,s)]=0
                 \end{equation*}
            \end{prop}
            \begin{prop}
                 \begin{equation*}
                      \{e_i(m,n),e_k(p,s)\}_+=\tau q^{m(s-n)}g_{ki}(m+p,s-n)
                 \end{equation*}
                 \begin{eqnarray*}
                      &&\{e_i(m,n),e_k^*(p,s)\}_+=\tau\delta_{ik}q^{-n(m+p)}e_0(m+p,s-n)+\tau q^{p(n-s)}f_{ik}(m+p,n-s)\\
                      &&\hspace{4cm}-\tau\delta_{ik}\delta_{m+p,0}\frac{1}{2}(q^{s-n}+1)\frac{q^{m(s-n)}-1}{q^{s-n}-1}
                 \end{eqnarray*}
                 \begin{equation*}
                      [e_i(m,n),e_0(p,s)]=-q^{np}e_i(m+p,n+s)+q^{p(n-s)}e_i(m+p,n-s)
                 \end{equation*}
            \end{prop}
            \begin{prop}
                 \begin{equation*}
                      \{e_i^*(m,n),e_k^*(p,s)\}_+=\tau q^{m(s-n)}h_{ki}(m+p,s-n)
                 \end{equation*}
                 \begin{equation*}
                      [e_i^*(m,n),e_0(p,s)]=-q^{np}e_i^*(m+p,n+s)+q^{p(n-s)}e_i^*(m+p,n-s)
                 \end{equation*}
            \end{prop}
            \begin{prop}
                 \begin{eqnarray*}
                      &&[e_0(m,n),e_0(p,s)]    \\
                      &=&-(q^{np}-q^{sm})e_0(m+p,n+s)
                          +\delta_{m+p,0}q^{np}\frac{1}{2}(q^{n+s}+1)\frac{q^{m(n+s)}-1}{q^{n+s}-1} \\
                      &&-(q^{m(s-n)}-q^{-n(m+p)})e_0(m+p,s-n)-
                         \delta_{m+p,0}\frac{1}{2}(q^{s-n}+1)\frac{q^{m(s-n)}-1}{q^{s-n}-1}
                 \end{eqnarray*}
            \end{prop}

       We only give proofs for 2.10 and 2.12. The
       proof for the others is either similar or easy.

       \noindent\textbf{Proofs of 2.10 and 2.12:} \\

Note that from [CG], we have
            \begin{equation}
            \sum_{t \in\mathbb{Z}}q^{-xt}\Bigl(\theta(-2t)-\theta(-2m-2t)\Bigr)
            =\frac{1}{2}(q^x+1)\frac{q^{mx}-1}{q^x-1}.
            \end{equation}

By (2.30),(2.31),(2.32) and (2.35),  we have
            \begin{equation}
                 \{e_i(m,n),a_k(p)e(s)\}_+=\tau q^{ns}a_k(p)a_i(m+s),
            \end{equation}
            \begin{equation}
                 \{e_i(m,n),a_k^*(p)e(s)\}_+=\tau q^{ns}a_k^*(p)a_i(m+s)+\tau \delta_{ik}q^{-n(m+p)}e(m+p)e(s),
            \end{equation}
            \begin{equation}
                 [e_i(m,n),e(p)e(s)]=q^{ns}a_i(m+s)e(p)-q^{np}a_i(m+p)e(s),
            \end{equation}
       and
            \begin{equation}
                 [e_0(m,n),e(p)e(s)]=-(q^{np}-q^{-n(m+p)})e(m+p)e(s)-(q^{ns}-q^{-n(m+s)})e(p)e(m+s),
            \end{equation}
       so we get
            \begin{eqnarray*}
                 &&\{e_i(m,n),q^{-st}:a_k(p-t)e(t):\}_+=\tau q^{nt-st}a_k(p-t)a_i(m+t)        \\
                 &=&\tau q^{m(s-n)}q^{-(s-n)(m+t)}a_k(p-t)a_i(m+t)                  \\
                 &=&\tau
                 q^{m(s-n)}q^{-(s-n)(m+t)}:a_k(p-t)a_i(m+t):,
            \end{eqnarray*}
            \begin{eqnarray*}
                 &&\{e_i(m,n),q^{-st}:a_k^*(p-t)e(t):\}_+                                \\
                 &=&q^{-st}\Bigl(\tau \delta_{ik}q^{-n(m+p-t)}e(m+p-t)e(t)+\tau q^{nt}a_k^*(p-t)a_i(m+t)\Bigr)     \\
                 &=&\tau\delta_{ik}q^{-n(m+p)}q^{-(s-n)t}e(m+p-t)e(t)+\tau q^{-(s-n)t}a_k^*(p-t)a_i(m+t)            \\
                 &=&\tau\delta_{ik}q^{-n(m+p)}q^{-(s-n)t}\Bigl(:e(m+p-t)e(t):-\delta_{m+p,0}\theta(m+p-2t)\Bigr)   \\
                 &&\tau q^{-(s-n)t}\bigl(:a_i(m+t)a_k^*(p-t):+\delta_{ik}\delta_{m+p,0}\theta(p-m-2t)\bigr)        \\
                 &=&\tau\delta_{ik}q^{-n(m+p)}q^{-(s-n)t}:e(m+p-t)e(t):
                            +\tau q^{p(n-s)}q^{-(n-s)(p-t)}:a_i(m+t)a_k^*(p-t):\\
                 &&-\tau\delta_{ik}\delta_{m+p,0}q^{-(s-n)t}\bigl(\theta(-2t)-
                 \theta(-2m-2t) \bigr);
            \end{eqnarray*}
            \begin{eqnarray*}
                 &&[e_i(m,n),q^{-st}:e(p-t)e(t):]      \\
                 &=&q^{-st}\bigl(q^{n(p-t)}a_i(m+p-t)e(t)+q^{nt}a_i(m+t)e(p-t)\bigr)      \\
                 &=&-q^{np}q^{-(n+s)t}:a_i(m+p-t)e(t):+q^{p(n-s)}q^{-(n-s)(p-t)}:a_i(m+t)e(p-t):     \\
            \end{eqnarray*}
       and
            \begin{eqnarray*}
                 &&[e_0(m,n),q^{-st}:e(p-t)e(t):]    \\
                 &=&-q^{-st}\Bigl((q^{n(p-t)}-q^{-n(m+p-t)})e(m+p-t)e(t)-(q^{nt}-q^{-n(m+t)}) e(p-t)e(m+t)\Bigr) \\
                 &=&-q^{-st}(q^{n(p-t)}-q^{-n(m+p-t)})\bigl(:e(m+p-t)e(t):-\delta_{m+p,0}\theta(m+p-2t)\bigr)\\
                 &&-q^{-st}(q^{nt}-q^{-n(m+t)})\bigl(:e(p-t)e(m+t):-\delta_{m+p,0}\theta(p-m-2t)\bigr) \\
                 &=&-q^{np}q^{-(s+n)t}:e(m+p-t)e(t):+q^{-n(m+p)}q^{-(s-n)t}:e(m+p-t)e(t):\\
                  &&-q^{m(s-n)}q^{-(s-t)(m+t)}:e(p-t)e(m+t):+q^{sm}q^{-(n+s)(m+t)}:e(p-t)e(m+t): \\
                  &&+\delta_{m+p,0}q^{np}q^{-(n+s)t}\bigl(\theta(-2t)-\theta(-2m-2t)\bigr)  \\
                  &&-\delta_{m+p,0}q^{-(s-n)t}\bigl(\theta(-2t)-\theta(-2m-2t)\bigr)
            \end{eqnarray*}
       by (2.36), we see that Proposition 2.10 and 2.12 hold true.
       $\hfill \blacksquare$

       To find the correspondence of the homomorphism, we need to
       modify our operators by shifting some central elements.

       For Proposition 2.3, we see that, if $n+s \in \Lambda(q)$
       and $n-s \in \Lambda(q)$,
            \begin{eqnarray*}
                  &&[g_{ij}(m,n),h_{kl}(p,s)]\\
                  &=&-\delta_{ik}q^{-n(m+p)}f_{jl}(m+p,s-n)-\delta_{jk}q^{np}f_{il}(m+p,n+s)  \\
                  &&-\delta_{il}q^{-(mn+np+ps)}f_{jk}(m+p,-(n+s))-\delta_{jl}q^{(n-s)p}f_{ik}(m+p,n-s) \\
                  &&+\delta_{ik}\delta_{jl}\delta_{m+p,0}m
                    +\delta_{jk}\delta_{il}\delta_{m+p,0}q^{np}m.
            \end{eqnarray*}

       If $n+s \in \mathbb{Z} \setminus \Lambda(q)$ and $n-s \in
       \Lambda(q)$, then
            \begin{eqnarray*}
                  &&[g_{ij}(m,n),h_{kl}(p,s)]\\
                  &=&-\delta_{ik}q^{-n(m+p)}f_{jl}(m+p,s-n)-\delta_{jl}q^{(n-s)p}f_{ik}(m+p,n-s)                 \\
                  &&-\delta_{jk}q^{np}\Bigl(f_{il}(m+p,n+s)
                       +\frac{1}{2}\delta_{il}\delta_{m+p,0}\frac{q^{n+s}+1}{q^{n+s}-1}\Bigr)   \\
                  &&-\delta_{il}q^{-(mn+np+ps)}\Bigl(f_{jk}(m+p,-n-s)
                       +\frac{1}{2}\delta_{jk}\delta_{m+p,0}\frac{q^{-n-s}+1}{q^{-n-s}-1})\Bigr)  \\
                  &&+\delta_{ik}\delta_{jl}\delta_{m+p,0}m.
            \end{eqnarray*}

       Similarly, if $n+s \in \Lambda(q)$ and $n-s \in \mathbb{Z}
       \setminus \Lambda(q)$
            \begin{eqnarray*}
                  &&[g_{ij}(m,n),h_{kl}(p,s)]      \\
                  &=&-\delta_{jk}q^{np}f_{il}(m+p,n+s)-\delta_{il}q^{-(mn+np+ps)}f_{jk}(m+p,-n-s)       \\
                  &&-\delta_{ik}q^{-n(m+p)}\Bigl(f_{jl}(m+p,s-n)
                       +\frac{1}{2}\delta_{jl}\delta_{m+p,0}\frac{q^{s-n}+1}{q^{s-n}-1}\Bigr)   \\
                  &&-\delta_{jl}q^{(n-s)p}\Bigl(f_{ik}(m+p,n-s)
                       +\frac{1}{2}\delta_{ik}\delta_{m+p,0}\frac{q^{n-s}+1}{q^{n-s}-1})\Bigr)  \\
                  &&+\delta_{jk}\delta_{il}\delta_{m+p,0}q^{np}m.
            \end{eqnarray*}

       By the above two relations, we have if $n+s,n-s \in
       \mathbb{Z} \setminus \Lambda(q)$
            \begin{eqnarray*}
                  &&[g_{ij}(m,n),h_{kl}(p,s)]      \\
                  &=&-\delta_{jk}q^{np}\Bigl(f_{il}(m+p,n+s)
                       +\frac{1}{2}\delta_{il}\delta_{m+p,0}\frac{q^{n+s}+1}{q^{n+s}-1}\Bigr)   \\
                  &&-\delta_{il}q^{-(mn+np+ps)}\Bigl(f_{jk}(m+p,-n-s)
                       +\frac{1}{2}\delta_{jk}\delta_{m+p,0}\frac{q^{-n-s}+1}{q^{-n-s}-1})\Bigr)  \\
                  &&-\delta_{ik}q^{-n(m+p)}\Bigl(f_{jl}(m+p,s-n)
                       +\frac{1}{2}\delta_{jl}\delta_{m+p,0}\frac{q^{s-n}+1}{q^{s-n}-1}\Bigr)   \\
                  &&-\delta_{jl}q^{(n-s)p}\Bigl(f_{ik}(m+p,n-s)
                       +\frac{1}{2}\delta_{ik}\delta_{m+p,0}\frac{q^{n-s}+1}{q^{n-s}-1})\Bigr)
            \end{eqnarray*}

       Using the same method, for Proposition 2.5 we have, if $n+s
       \in \Lambda(q)$,
            \begin{eqnarray*}
                 &&[f_{ij}(m,n),f_{kl}(p,s)]    \\
                 &=&\delta_{jk}q^{np}f_{il}(m+p,n+s)-\delta_{il}q^{sm}f_{kj}(m+p,n+s)
                       -\delta_{jk} \delta_{il} q^{np} \delta_{m+p,0} m .
            \end{eqnarray*}

       If $n+s \in \mathbb{Z} \setminus \Lambda(q)$, then
            \begin{eqnarray*}
                 &&[f_{ij}(m,n),f_{kl}(p,s)]    \\
                 &=&\delta_{jk}q^{np}\Bigl(f_{il}(m+p,n+s)
                      +\frac{1}{2}\delta_{il}\delta_{m+p,0}\frac{q^{n+s}+1}{q^{n+s}-1}\Bigr)   \\
                 &&-\delta_{il}q^{sm}\Bigl(f_{kj}(m+p,n+s)
                      +\frac{1}{2}\delta_{jk}\delta_{m+p,0}\frac{q^{n+s}+1}{q^{n+s}-1})\Bigr).
            \end{eqnarray*}

       For Proposition 2.10, if $n-s \in \Lambda(q)$,
            \begin{eqnarray*}
                 &&\{e_i(m,n),e_k^*(p,s)\}_+    \\
                 &=&\tau\delta_{ik}q^{-n(m+p)}e_0(m+p,s-n)+\tau q^{p(n-s)}f_{ik}(m+p,n-s)
                      -\tau\delta_{ik}\delta_{m+p,0}m .
            \end{eqnarray*}

       If $n+s \in \mathbb{Z} \setminus \Lambda(q)$, then
            \begin{eqnarray*}
                 &&\{e_i(m,n),e_k^*(p,s)\}_+    \\
                 &=&\tau\delta_{ik}q^{-n(m+p)}\Bigl(e_0(m+p,s-n)
                      +\frac{1}{2}\delta_{m+p,0}\frac{q^{s-n}+1}{q^{s-n}-1}\Bigr) \\
                 &&+\tau q^{p(n-s)}\Bigl(f_{ik}(m+p,n-s)
                      +\frac{1}{2}\delta_{jk}\delta_{m+p,0}\frac{q^{n-s}+1}{q^{n-s}-1})\Bigr).
            \end{eqnarray*}

       For Proposition 2.12, if $n+s \in \Lambda(q)$ and $n-s \in
       \Lambda(q)$,
            \begin{eqnarray*}
                 &&[e_0(m,n),e_0(p,s)]    \\
                 &=&-(q^{np}-q^{sm})e_0(m+p,n+s)
                       -(q^{m(s-n)}-q^{-n(m+p)})e_0(m+p,s-n) \\
                 &&+\delta_{m+p,0}q^{np}m-\delta_{m+p,0}m.
            \end{eqnarray*}

       If $n+s \in \mathbb{Z} \setminus \Lambda(q)$ and $n-s \in
       \Lambda(q)$, then
            \begin{eqnarray*}
                 &&[e_0(m,n),e_0(p,s)]    \\
                 &=&-(q^{np}-q^{sm})\Bigl(e_0(m+p,n+s)
                      +\frac{1}{2}\delta_{m+p,0}\frac{q^{n+s}+1}{q^{n+s}-1}\Bigr) \\
                 &&-(q^{m(s-n)}-q^{-n(m+p)})e_0(m+p,s-n)-\delta_{m+p,0}m.
            \end{eqnarray*}

       If $n+s \in \Lambda(q)$ and $n-s \in \mathbb{Z} \setminus
       \Lambda(q)$
            \begin{eqnarray*}
                 &&[e_0(m,n),e_0(p,s)]    \\
                 &=&-(q^{np}-q^{sm})e_0(m+p,n+s)+\delta_{m+p,0}q^{np}m                \\
                 &&-(q^{m(s-n)}-q^{-n(m+p)})\Bigl(e_0(m+p,s-n)
                      +\frac{1}{2}\delta_{m+p,0}\frac{q^{s-n}+1}{q^{s-n}-1}\Bigr).
            \end{eqnarray*}

       If $n+s,n-s \in \mathbb{Z} \setminus \Lambda(q)$
            \begin{eqnarray*}
                 &&[e_0(m,n),e_0(p,s)]    \\
                 &=&-(q^{np}-q^{sm})\Bigl(e_0(m+p,n+s)
                      +\frac{1}{2}\delta_{m+p,0}\frac{q^{n+s}+1}{q^{n+s}-1}\Bigr)           \\
                 &&-(q^{m(s-n)}-q^{-n(m+p)})\Bigl(e_0(m+p,s-n)
                      +\frac{1}{2}\delta_{m+p,0}\frac{q^{s-n}+1}{q^{s-n}-1}\Bigr).
            \end{eqnarray*}

       Therefore, if we define \\
            \noindent \parbox{1cm}{\begin{eqnarray}\end{eqnarray}}\hfill \parbox{13.66cm}
               {\begin{eqnarray*}& F_{ij}(m,n)=\left\{
                   \begin{array}
                       {r@{ \quad }l}  f_{ij}(m,n), & \mbox{  for } n \in \Lambda(q)  \\
                       f_{ij}(m,n)+\frac{1}{2}\delta_{ij}\delta_{m,0}\frac{q^n+1}{q^n-1},
                       &\mbox{  for } n \in \mathbb{Z} \setminus \Lambda(q)
                   \end{array} \right. & \\
                       &G_{ij}(m,n)=g_{ij}(m,n),\quad H_{ij}(m,n)=h_{ij}(m,n),& \\
                       &E_i(m,n)=e_i(m,n),\quad E_i^*(m,n)=e_i^*(m,n),& \\
                       & E_0(m,n)=\left\{
                   \begin{array}
                        {r@{ \quad }l} e_0(m,n), & \mbox{  for } n \in \Lambda(q) \vspace{0.1cm} \\
                        e_0(m,n)+\frac{1}{2}\delta_{m,0}\frac{q^n+1}{q^n-1},
                        &\mbox{  for } n \in \mathbb{Z} \setminus \Lambda(q)
                   \end{array} \right. & \\ \vspace{1cm}  \hfill
                \end{eqnarray*}}
                we have
            \begin{theorem}
                $V_\tau(N)$ is a module for the Lie superalgebra
                $\widehat{\mathcal{G}}$ under the action given by
                (for $\tau=\pm1$)
                \begin{eqnarray*}
                     \pi(\tilde{g}_{ij}(m,n))=\tau G_{ij}(m,n), &&\qquad  \pi(\tilde{f}_{ij}(m,n))=F_{ij}(m,n),\\
                     \pi(\tilde{h}_{ij}(m,n))=\tau H_{ij}(m,n), &&\qquad  \pi(\tilde{e}_i(m,n))=\tau E_i(m,n),   \\
                     \pi(\tilde{e}_i^*(m,n))=E_i^*(m,n),        &&\qquad  \pi(\tilde{e}_0(m,n))=E_0(m,n),      \\
                     \pi(c_x)=-\frac{1}{2},                   &&\qquad \pi(c_y)=0.
                \end{eqnarray*}
            \end{theorem}

 \footnotesize{    \begin{align*}
        &\mbox{Hongjia Chen, Shikui Shang}                                    &&\mbox{Yun Gao}          \\
         &\mbox{Department of Mathematics}                       &&\mbox{Department of Mathematics and Statistics} \\
         &\mbox{University of Science and Technology of China}   &&\mbox{York University} \\
         &\mbox{Hefei, Anhui}                                    &&\mbox{Toronto, Ontario} \\
         &\mbox{P. R. China  230026}                             &&\mbox{Canada  M3J 1P3} \\
         &\mbox{Email:hjchen@mail.ustc.edu.cn}                   &&\mbox{Email:ygao@yorku.ca} \\
         &\mbox{Email:skshang@mail.ustc.edu.cn}
       \end{align*}}
\end{document}